\documentclass[11pt]{article}
\usepackage{amsfonts,amsmath,amssymb,amsthm, amscd}
\usepackage{graphicx}

\numberwithin{equation}{section}

\newtheorem{theorem}{Theorem}[section]

\theoremstyle{definition}

\newtheorem{remark}[theorem]{Remark}

\def\<{{\langle}}
\def\>{{\rangle}}
\def\G{{\Gamma}}
\def\a{{\alpha}}
\def\b{{\beta}}
\def\g{{\gamma}}
\def\Z{\mathbb Z}

\def\R{\mathbb R}

\def\S{{\mathbb S}}

\def\a{\alpha}

\def\G{{\Gamma}}

\def\ni{\noindent} 

\begin{document}

\title{On the Component Number of Links from Plane Graphs}

\author{Daniel S. Silver 
\and Susan G. Williams\thanks {Both authors are partially supported by the Simons Foundation.} }

\maketitle 

\begin{abstract} A short, elementary proof is given of the result: The number of components of a link arising from a medial graph $M(\G)$ by resolving vertices is equal to the nullity of the mod-$2$ Laplacian matrix of $\G$. 
\end{abstract}


\section{Introduction} \label{Intro} 
\begin{figure}
\begin{center}
\includegraphics[height=1.5 in]{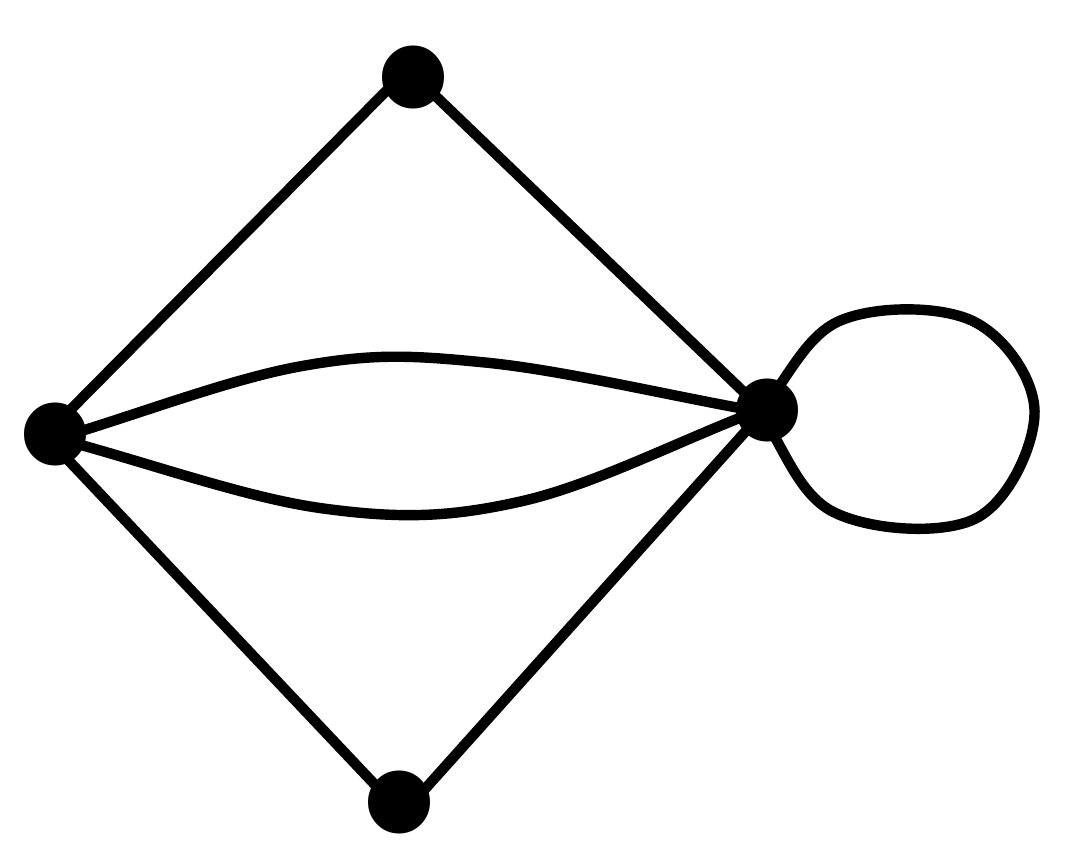}
\caption{Plane Graph $\G$}
\label{graph}
\end{center}
\end{figure}

\begin{figure}
\begin{center}
\includegraphics[height=1.5 in]{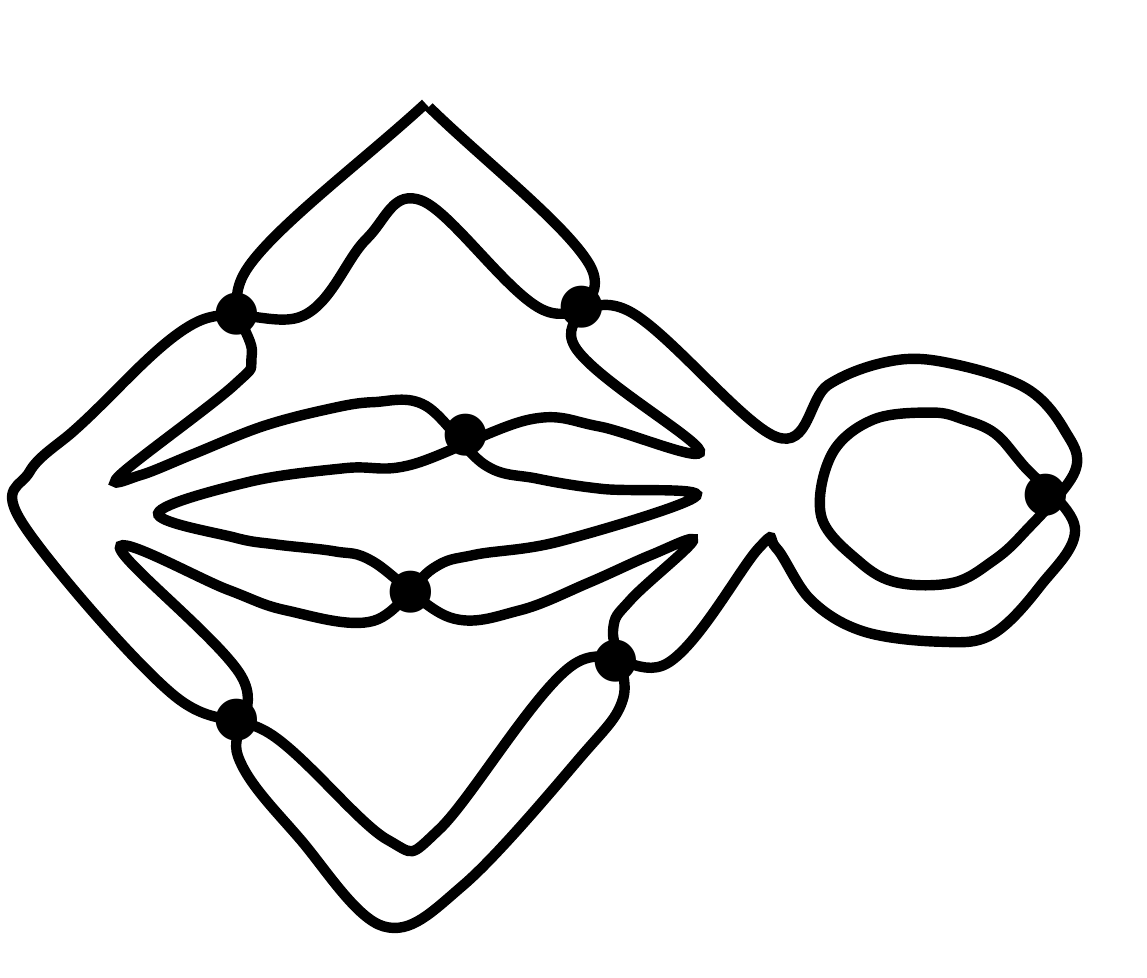}
\caption{Medial graph $M(\G)$}
\label{medial}
\end{center}
\end{figure}

Let $\G$ be a plane graph. Its \emph{medial graph} $M(\G)$ is obtained from the boundary of a regular neighborhood of $\G$ by pinching  each edge to create a vertex of degree 4. An example is given in Figures \ref{graph} and \ref{medial}. 

\begin{figure}
\begin{center}
\includegraphics[height=1.5 in]{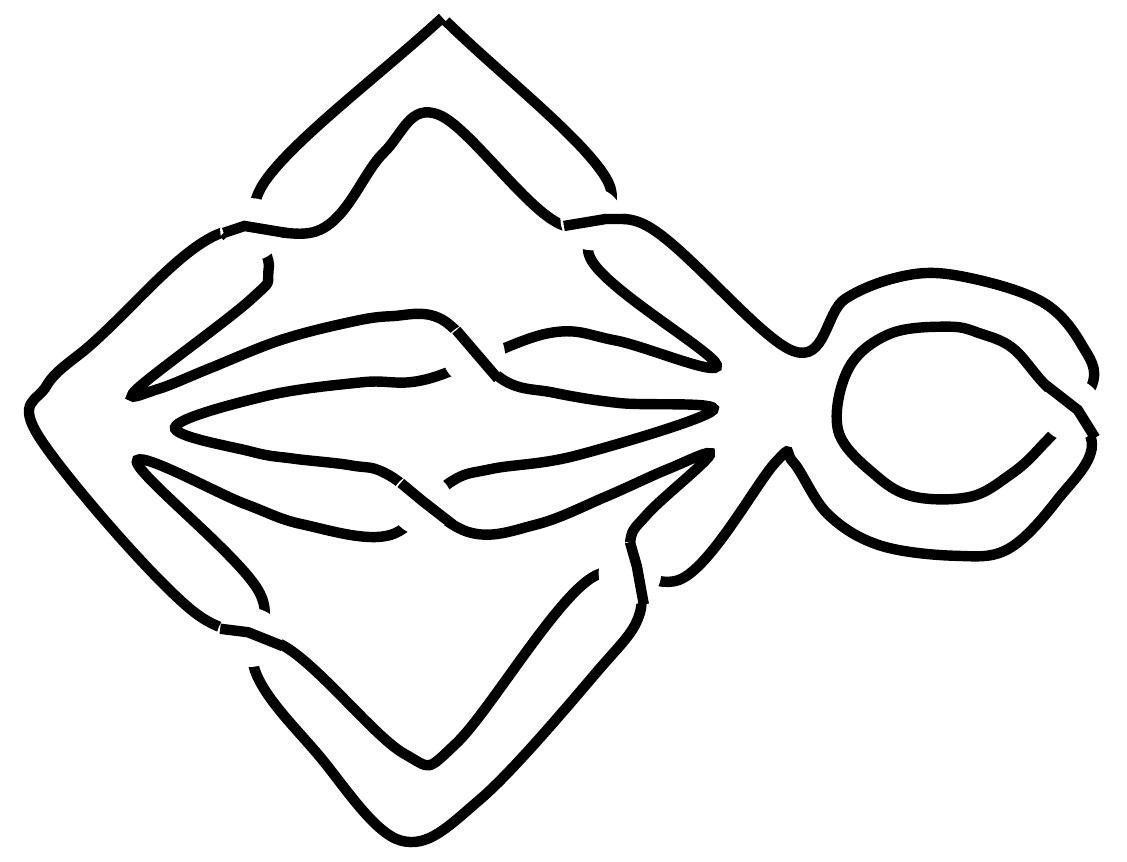}
\caption{Link $L$ associated to $M(\G)$}
\label{link}
\end{center}
\end{figure}

The construction is important in knot theory, since any diagram of a link can obtained from a suitable plane graph $\G$ by resolving each vertex of $M(\G)$ in one of two ways so that one arc of the diagram appears to pass over the other.  (See Figure \ref{link}.) In this way, questions about knots and links can be converted into the language of graph theory. Indeed, many invariants of the link can be computed directly from $\G$ (see, for example, \cite{kauffman}, \cite{kauffman2}). 

The most basic invariant of a link $L$ is the number of its components, denoted here by $\mu(L)$. Notice that $\mu(L)$ is the same for all links associated to a medial graph $M(\G)$, regardless of how the vertices are resolved. Calculation of $\mu(L)$ for links arising from special families of graphs have appeared in various places (for example, \cite{mphako}, \cite{ptz}).

In \cite{lv} it is shown that the Tutte polynomial evaluation $T_\G(-1,-1)$  determines $\mu(L)$. A more intrinsic formula appears in \cite{gr}: 

The \emph{mod-2 Laplacian} matrix of a graph $\G$ is the matrix $Q_2(\G)=(q_{ij})$, where $q_{ii}$ is the degree of the $i$th vertex, while for $i\ne j$,  $q_{ij}$ is the number of edges between the $i$th and $j$th vertices. All entries are taken modulo 2.  The following theorem is proved in \cite{gr}. 

\begin{theorem}\label{main} Let $L$ be a link arising from a medial graph $M(\G)$ by resolving vertices. The number $\mu(L)$ of components of $L$ is the nullity of the mod-2 Laplacian matrix $Q_2(\G)$. \end{theorem}

\begin{remark} As observed in \cite{lv}, both the Tutte polynomial and the Laplacian matrix determination of $\mu(\G)$ show that $\mu(L)$ is independent of the planar embedding of $\G$. \end{remark} 

Theorem \ref{main} is part of topology's folklore. It is provable by indirect, non-elementary techniques of algebraic topology.  One begins with the observation that $Q_2(\G)$ is the ``unreduced mod-$2$ Goertiz matrix" of the link $L$, a matrix shown by H. Seifert \cite{se} and R.H. Kyle \cite{ky}  to be a presentation matrix for
$H_1(M_2; \Z_2)\oplus \Z_2$, where $M_2$ is the 2-fold cyclic cover of $\S^3$ branched over $L$ (see, for example, \cite{lickorish}). The rank of $H_1(M_2; \Z_2)$ is known to be $\mu(L)-1$ (see, for example, \cite{gl}). Hence the nullity of $Q_2(\G)$ is $\mu(L)$. 
The complete argument can be found in \cite{na}.  

A combinatorial proof of Theorem \ref{main} is given in Chapter 17 of \cite{gr}. It draws on substantial material about cut, flow and bicycle spaces contained in earlier chapters.  A proof given in \cite{bkr} is more direct, making explicit the connection between the mod-$2$ Laplacian matrix of $Q_2(G)$ and the components of the medial graph $M(G)$. 

Here we add an idea of conservative vertex-coloring, motivated by knot theory, to give an extremely short and self-contained proof of Theorem \ref{main}. The idea has been adapted for locally-finite periodic graphs in \cite{lsw}. 

The authors wish to thank Iain Moffatt and Lorenzo Traldi for helpful suggestions. 

\section{Flats}   
\begin{figure}
\begin{center}
\includegraphics[height=1.5 in]{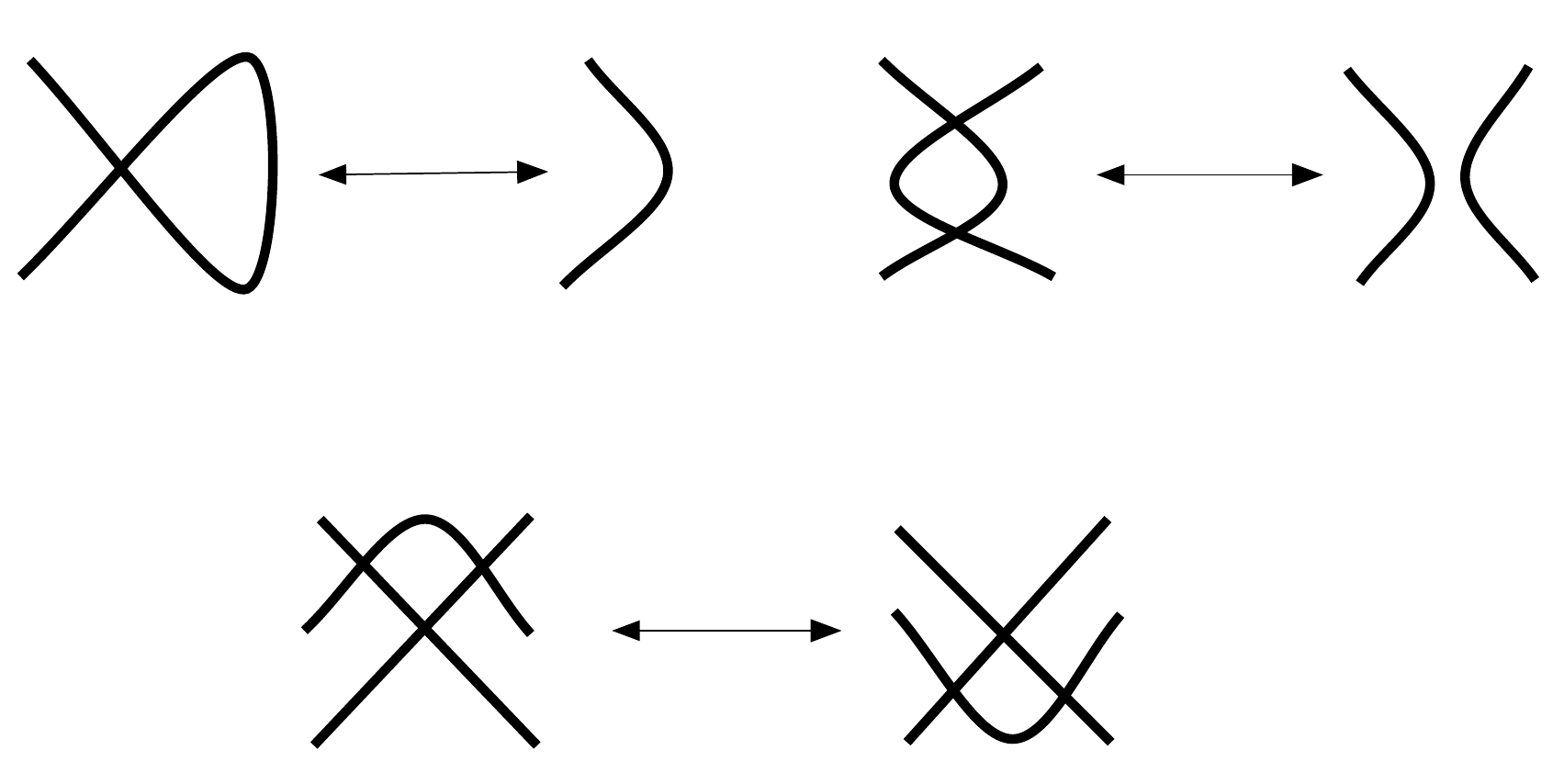}
\caption{Flat Reidemeister moves}
\label{flat}
\end{center}
\end{figure}

A \emph{flat} $F$ is a finite collection of immersed circles in the plane that are in \emph{general position} (that is,   at any intersection point, the two tangent vectors span the plane).  Two flats are \emph{equivalent} if one can be changed into the other by a finite sequence of \emph{flat Reidemeister moves}, shown in Figure \ref{flat}. Such moves are understood to performed locally, fixing any points not shown. We also allow homeomorphisms of the plane, which permit deformations such as bending and stretching. 

Flats are produced by link diagrams by projection. Not surprisingly, flat Reidemeister moves were inspired by well-known Reidemeister moves of link diagrams in the plane. (Both sets of moves are described in \cite{kauffman}).  We leave it as the first of two elementary exercises to show that any flat  $F$ is equivalent to a collection of circles that bound pairwise disjoint embedded disks in the plane. We call the number of  circles the \emph{component number} of $F$, and we denote it by $\mu(F)$. It is clearly well defined.

Given a flat $F$, we associate a vector space $V_F$ over the field $\Z_2$ by means of a presentation. 
The generators correspond to the bounded regions of $\R^2 \setminus F$. At each crossing we impose the relation that the sum of the generators corresponding to the regions meeting there is zero.  

For the second exercise, one verifies that if two flats differ by a flat Reidemeister move, then the corresponding vector spaces are isomorphic. Consequently, $V_F$ has dimension $\mu(F)$.


\section{Proof of Theorem \ref{main}} Let $\G$ be a plane graph, $M(\G)$ its medial graph. By resolving vertices, various diagrams of links are possible. They all have $\mu_F$ components, where $\mu_F$ is the component number of the underlying flat $F$.

Checkerboard shade the regions of $\R^2\setminus F$ with the unbounded region unshaded. Then $\G$ can be identified with the graph having a vertex in each shaded region and an edge running through each crossing.  

Recall that a vector space $V_F$ is associated to $F$ by means of a presentation ${\cal P}_F$. By \emph{shaded} (resp. \emph{unshaded})  \emph{generators} of ${\cal P}_F$ we mean  generators corresponding to shaded (resp. unshaded) regions of $\R^2\setminus F$. We identify shaded generators with vertices of $\G$, and the unshaded generators with faces of $\G$. Since the unbounded region corresponds to the trivial element, it is easy to see from the relations that the shaded generators alone generate $V_F$. 

A \emph{vertex coloring} is an assignment of colors in $\Z_2$ to the shaded generators. It is \emph{conservative} if it extends over unshaded regions to give an element of $V_F$.  If a vertex coloring extends, then it does so uniquely. Hence we can count the elements of $V_F$ by counting conservative vertex colorings.

Assume that $v$ is a vertex of $\G$ with valence $d$. Let $v_1, \ldots, v_d$ be the vertices adjacent to $v$. (We allow repetition in the case of multiple edges.) Let $\a_v$ (resp. $\a_{v_i}$) be  colors assigned to to $v$ (resp. $v_i$). One checks that if the vertex coloring is conservative, then the $d$ relations $v$ imply: 
\begin{equation}\label{first} \a_{v_1} + \ldots + \a_{v_d} =0,\quad (d\ {\rm even})
\end{equation}
\begin{equation}\label{second} \a_v + \a_{v_1} + \ldots + \a_{v_d} =0,\quad (d\ {\rm odd}) \end{equation}

Conversely, if the above conditions hold at every vertex, then the coloring is a conservative vertex coloring. To prove this, it suffices to show that the coloring extends to the unshaded regions in a way that is consistent with all relations. For this, consider a path from the unbounded region to any unshaded region $R$, intersecting edges of $\G$ transversely. Determine colors for the unshaded regions inductively: as the path leaves a region labeled $\g$ and crosses an edge of $\G$ with vertices labeled $\a$ and $\b$,  assign $\a + \b + \g$ to the unshaded region that the path enters. Then the  relation holds at the crossing corresponding to that edge. We call $\a + \b + \g$ the result of \emph{integrating along the path}. The conditions (\ref {first}) and (\ref {second}) imply that if we integrate around a small {\sl closed} path surrounding any vertex, then initial and final colors agree. Consequently, integrating along any two paths from the unbounded region to $R$ yields the same result. Since the path can be taken to cross any given edge, all  relations are satisfied.

In view of the equations that must hold at each vertex, conservative vertex colorings are the elements of the null space of the mod-2 Laplacian matrix $Q_2(\G)$. Hence Theorem \ref{main} is proved. 

\begin{remark} If $L$ is any link associated to $M(\G)$ by resolving vertices, then the \emph{Dehn presentation} of $\pi_1(\R^3\setminus L)$ has generators corresponding to the bounded regions of $\R^2 \setminus D$, where $D$ is the diagram of $L$, and relations as 
in Figure \ref{Dehn}. Our proof is inspired by the observation that 
$V_L$ can be identified with ${\rm Hom}(\pi_1(\R^3 \setminus L), \Z_2)$. 

\begin{figure}
\begin{center}
\includegraphics[height=1.5 in]{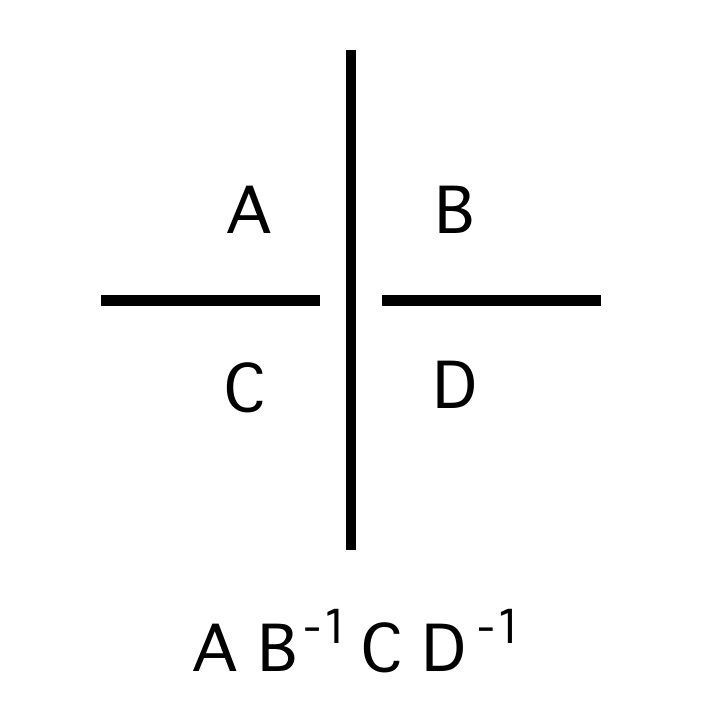}
\caption{Dehn generators and relation}
\label{Dehn}
\end{center}
\end{figure}

\end{remark}


\bigskip

\ni Department of Mathematics and Statistics,\\
\ni University of South Alabama\\ Mobile, AL 36688 USA\\
\ni  silver@southalabama.edu\\
\ni swilliam@southalabama.edu
\end{document}